\documentclass[12pt]{article}

\usepackage[english]{babel}
\usepackage{graphicx}
\usepackage{mathtools}
\usepackage[toc,page]{appendix}
\usepackage{amsmath,amsfonts,latexsym }%,showkeys}

\usepackage{pdfpages}
\usepackage{amssymb}
\usepackage{graphics}
\usepackage{geometry}
\usepackage{lipsum}     
\usepackage[toc,page]{appendix}
\usepackage{fancyhdr}
\usepackage{setspace}

\newcommand{\nc}{\newcommand}
\nc{\tr}{{\vartriangle}} \nc{\vth}{{\vartheta}}
\nc{\bt}{{\beta}}
\nc{\dl}{{\delta}} \nc{\Dl}{{\tr}}
\nc{\p}{{\psi}}
\nc{\gm}{{\gamma}} \nc{\Gm}{{\Gamma}} \nc{\sg}{{\sigma}}
\nc{\ve}{{\varepsilon}} 
\nc{\ch}{{\cal H}} 
\nc{\cf}{{\cal F}}
\nc{\cp}{{\cal P}}
\nc{\td}{\tilde}
\nc{\ck}{\cal K}

\newtheorem{lemma}{Lemma}[section]
\newtheorem{theorem}[lemma]{Theorem}
\newtheorem{corollary}[lemma]{Corollary}
\newtheorem{proposition}[lemma]{Proposition}

\newtheorem{example}[lemma]{Example}

\newtheorem{remark}[lemma]{Remark}
\numberwithin{equation}{section}

\begin{document}

\title
{Asymptotic Behaviour of Truncated Stochastic Approximation Procedures}
   \author{Teo Sharia and Lei Zhong}
\date{}
\maketitle
\begin{center}
{\it \footnotesize
Department of Mathematics, Royal Holloway,  University of London\\
Egham, Surrey TW20 0EX \\ e-mail: t.sharia@rhul.ac.uk }
\end{center}

%\thankstext{T1}{Footnote to the title with the `thankstext' command.}

\begin{abstract}
We study asymptotic behaviour of  stochastic approximation procedures
with   three main characteristics:  truncations  with random moving bounds, 
 a matrix valued random step-size sequence,  and  a dynamically changing random regression function.
 In particular, we  show   that under quite mild conditions, stochastic approximation procedures are asymptotically linear in the statistical sense, that is, they can be represented as weighted sums of random variables. Therefore, a suitable form of the central limit theorem can be applied to derive asymptotic distribution of the corresponding processes.  The  theory  is illustrated by various examples and special cases.
  \end{abstract}

%\begin{keyword}[class=AMS]
%\kwd[Primary ]{{62F10, 62F12, 62M10}}
%\kwd[; secondary ]{62L20, 62F35}
%\end{keyword}

\begin{center}
\end{center}

\section{Introduction}\label{Intro}
This paper  is the final part  of the series of papers devoted to the study of  truncated Stochastic approximation (SA)   with moving bounds.
The   classical problem  of SA is concerned with finding a unique zero, say $z^0$,  of  a real valued function $R(z):  \mathbb{R}  \to \mathbb{R}$ when only noisy measurements  of $R$ are available. To estimate $z^0$, consider a sequence defined recursively   as                         
$$
Z_t=
Z_{t-1}+\gm_t \left[R(Z_{t-1})+\ve_t\right],   \qquad t=1,2,\dots
$$
where $\{\ve_t\}$ is a sequence of zero-mean random variables and  $\{\gamma_t\}$ is a deterministic sequence of positive numbers.
This is the classical Robbins-Monro SA procedure (see  Robbins and Monro (1951)\nocite{RM}), which under certain conditions  converges to the root $z^0$ of the equation $R(z)=0$.  (Comprehensive surveys of the SA technique can be found in Benveniste et al. (1990)\nocite{Ben1990}, Borkar (2008)\nocite{Bor},
Kushner and Yin (2003)\nocite{KushYin}, Lai (2003)\nocite{Lai}, and Kushner (2010)\nocite{Kush1}.)  

In applications however, it is important to consider the setting  when the function $R$ changes over the time. So, let us assume that 
the objective now is to find a common root $z^0$ of a dynamically changing sequence of functions $R_t(z)$.  
Also, in certain circumstances it might be necessary to confine the values of the procedure to a certain set, or to a sequence of sets by applying a truncation operator.  This happens if,  e.g., the functions in the recursive equation are defined only for certain values of the parameter. Truncations  may  also be useful when certain  standard assumptions, e.g., 
  conditions on  the  growth rate of the relevant functions are not satisfied.   
Truncations may  also help to make an efficient use of   auxiliary information  concerning the value of the unknown parameter.
For example, we might have auxiliary information about the root $z^0$, e.g. a set, possibly time dependent, that  contains the value of the unknown root. 
In order to study these procedures in an  unified  manner,  we consider a SA of the following form 
$$                         
Z_t=\Phi_{U_t} \Big(~
Z_{t-1}+\gm_t(Z_{t-1}) \big[ R_t(Z_{t-1})+\ve_t(Z_{t-1}) \big]\Big),   \quad t=1,2,\dots
$$
where  $Z_0 \in \mathbb{R}^m$  is some starting value,  $R_t(z)$ is  a predictable process with the property that $R_t(z^0)=0$ for all $t$'s, $\gm_t(z)$ is a matrix-valued predictable step-size sequence,  $U_t   \subset \mathbb{R}^m$  is  a random sequence of truncation sets,   and  $\Phi$ is the truncation operator  which returns the procedure to  $U_t$ every time the updated value leaves the truncation set  (see Section \ref{MON} for details). 
These SA procedures have the  following main characteristics: (1) inhomogeneous  random functions $R_t$; (2) state dependent matrix valued random step-sizes;  (3)  truncations with random and  moving 
(shrinking or expanding) bounds.
The main motivation for these comes from parametric statistical 
applications: (1) is needed for recursive parameter estimation procedures for  non i.i.d. models;
(2) is required to guarantee asymptotic optimality and efficiency of  statistical estimation; (3) is needed for
various different adaptive truncations, in particular, for the ones arising by auxiliary estimators  (see  Sharia (2014)\nocite{Shar4}
for a more detailed discussions  of these extensions).

Note that the idea of truncations  goes back to Khasʹminskii and Nevelson (1972)\nocite{Khas} and  Fabian (1978)\nocite{Fab}  (see also 
 Chen and Zhu (1986)\nocite{Chen2}, Chen et al.(1987)\nocite{Chen1},  Andrad{\'o}ttir (1995)\nocite{Andr}, Sharia (1997)\nocite{Shar0},  Tadic (1997,1998)\nocite{Tadic1}\nocite{Tadic2}, Lelong (2008)\nocite{Lel}.
A comprehensive bibliography and some comparisons can be found in Sharia (2014)\nocite{Shar4}).

 Convergence of the above class of procedures was studied  in Sharia (2014) and the results on  rate of convergence
 were established in  
 Sharia and Zhong (2016).  In this paper, we derive further asymptotic properties  of these procedures.  
 In particular, we show that  under quite mild conditions, SA procedures are asymptotically linear in the statistical sense, that is, they can be represented as weighted sums of random variables. Therefore, a suitable form of the central limit theorem can be applied to derive asymptotic distribution of the corresponding SA process.  Since some of the conditions in the main statements might be  difficult to interpret,  we present   explanatory remarks and corollaries. 
We also discuss  the case of the classical SA and  demonstrate that  truncations with moving bounds make it possible  to use SA  even when the standard conditions on the function $R$ do not hold.
 Finally,
applications of the above results are discussed and some simulations are presented to illustrate the theoretical results of the paper.    Proofs of some technical parts are postponed to Appendices.

%%%%%%%%%%%%%%%%%%%%%%%%%%%%%%%%%%%%%%%%%%%%%%%%%%%%%%%%%%%%%%%%%%%%%%%%%%%%%%%

%%%%%%%%%%%%%%%%%%%%%%%%%%%%%%%%%%%%%%%%%%%%%%%%%%%%%%%%%%%%%%%%%%%%%%%%%%%%%%%%%

%%%%%%%%%%%%%%%%%%%%%%%%%%%%%%%%%%%%%%%%%%%%%%%%%%%%%%%%%%%%%%%%%%%%%%%%%%
                 %       SUB S E C T I O N    Main   Results
%%%%%%%%%%%%%%%%%%%%%%%%%%%%%%%%%%%%%%%%%%%%%%%%%%%%%%%%%%%%%%%%%%%%%%%%%%

 \section{Main results}\label{MR}

\subsection{Notation and preliminaries}\label{MON}

Let $(\Omega, ~ \cf,F=(\cf_t)_{t\geq 0}, ~P)$ be a stochastic basis satisfying the usual conditions. Suppose that 
for each $t=1,2, \dots$, we have 
$( {\cal{B}}  ( \mathbb{R}^m)  \times \cf )$-measurable functions
$$
\begin{array}{cl}
 R_t(z)= R_t(z,\omega) &:\mathbb{R}^m \times  \Omega   \to    \mathbb{R}^m  \\
 \ve_t(z)=\ve_t(z,\omega) &:\mathbb{R}^m \times  \Omega   \to    \mathbb{R}^m      \\
  \gamma_t(z)=\gamma_t(z,\omega)&:\mathbb{R}^m \times  \Omega   \to    \mathbb{R}^{m\times m}    
\end{array}
$$
such that 
 for each $z\in  \mathbb{R}^m$,  the   processes    $R_t(z) $  and  $\gamma_t(z)$ are  predictable, i.e.,
 $R_t(z) $  and  $\gamma_t(z)$ are $\cf_{t-1}$ measurable for each $t$.
 Suppose also that  
 for each $z\in  \mathbb{R}^m$,  the process $\ve_t(z) $   is a martingale difference, i.e., $\ve_t(z) $ 
 is  $\cf_{t}$ measurable and  $E\left\{\ve_t(z)\mid{\cal{F}}_{t-1}\right\}=0$.
 We also assume that 
$$
R_t(z^0)=0
$$  
for each $ t=1, 2, \dots $, where  $z^0 \in   \mathbb{R}^m$   is  a non-random vector.

Suppose that $h=h(z)$ is a real valued function  of
 $ z \in {{\mathbb{R}}}^m$.  Denote by $ h'(z)$  the row-vector
 of partial derivatives
of $h$ with respect to the components of $z$, that
is,
 $
  h'(z)=\left(\frac{{\partial}}{{\partial} z_1} h(z), \dots,
 \frac{{\partial}}{{\partial} z_m} h(z)\right).
 $
Also, we denote by  $h''(z)$ the  matrix of second partial derivatives.
 The $m\times m$ identity matrix is denoted by ${{\bf I}}$.
%(here $T$ means transposition). 
Denote by $[a]^+$ and $[a]^-$ the positive and negative parts of $a\in \mathbb R$, i.e. $[a]^+=\max(a,0)$ and $[a]^-=\min(a,0)$.

Let    $U \subset    \mathbb{R}^m$ is a closed convex set  and define a truncation  operator    as a function
$\Phi_U(z) : \mathbb{R}^m  \longrightarrow \mathbb{R}^m$, such that
$$
\Phi_U(z) =\begin{cases}
    z & \text{if} \;\; z\in U \\
     z^* & \text{if} \;\; z\notin U,
\end{cases}
$$
where $z^*$ is a point in 
$U$,  that minimizes   the distance  to $z$.

Suppose that   $z^0 \in   \mathbb{R}^m$.  We say that a random sequence   of sets $U_t =U_t(\omega)$ 
($ t=1,2, \dots $)  from $\mathbb{R}^m $  is {\underline{\bf admissible}}  for  $z^0$  if
\medskip

\noindent
$\bullet$ 
for each $t$ and  $\omega,$    $U_t(\omega)$   is a   closed convex  subset  of $ \mathbb{R}^m$;
 \\
$\bullet$ 
 for each   $t$ and $z \in   \mathbb{R}^m$, the truncation $\Phi_{U_t}(z) $  is  $ {\cal{F}}_{t}$ measurable; 
 \\
$\bullet$   $z^0\in U_t$ eventually, i.e.,  for almost all  $\omega$ there exist  $t_0(\omega)<\infty$
 such that $z^0\in U_t(\omega)$ whenever $t >t_0(\omega)$.

\medskip
Assume that $Z_0 \in \mathbb{R}^m$  is some starting value and  consider the procedure 
                          %(TSA)
\begin{equation}\label{TSA}
Z_t=
\Phi_{U_t}
\Big(
Z_{t-1}+\gm_t(Z_{t-1}) \Psi_t(Z_{t-1})\Big),   \quad t=1,2,\dots
\end{equation}
where  $U_t $   is {admissible}   for  $z^0$,
$$
\Psi_t(z)=R_t(z)+\ve_t(z),
$$
and   $R_t(z) $, $\ve_t(z)$, $\gm_t(z)$ are random fields  defined above. Everywhere in this work, we assume that
                       %(GTSA1)
\begin{equation}\label{GTSA1}
 E\left\{\Psi_t(Z_{t-1})\mid{\cal{F}}_{t-1}\right\}=R_t(Z_{t-1})
  \end{equation}
                     and
\begin{equation}\label{GTSA2}
 E\left\{\ve_t^T(Z_{t-1})\ve_t(Z_{t-1})\mid{\cal{F}}_{t-1}\right\}= 
 \left[E\left\{\ve_t^T(z)\ve_t(z)\mid{\cal{F}}_{t-1}\right\} \right] _{z=Z_{t-1}},
 \end{equation}
 and the conditional expectations   \eqref{GTSA1} and \eqref{GTSA2}  are assumed to be finite.
 \medskip

%
  %%%%%%%%%%%%%%%%%%%%%%%%%%%%%%%%%%%%%%%%%%%%%%%%%%%%%%%%%%%%%%%%%%%%%%%%%%
                
%%%%%%%%%%%%%%%%%%%%%%%%%%%%%%%%%%%%%%%%%%%%%%%%%%%%%%%%%%%%%%%%%%%%%%%%%%
\begin{remark}\label{disint} {\rm  Condition \eqref{GTSA1} ensures that $\ve_t(Z_{t-1})$ is a martingale difference.
Conditions  \eqref{GTSA1} and \eqref{GTSA2} obviously hold if, e.g.,  the  measurement errors $\ve_t(u)$ are independent 
random variables, or if they  are state independent. In general,  
since we assume that all  conditional  expectations are calculated  as integrals w.r.t. corresponding regular conditional probability measures (see the convention below), these conditions can be checked using disintegration formula (see, e.g.,  Theorem 5.4 in Kallenberg (2002)\nocite{Kall}).}
\end{remark}
\noindent
{\bf \em Convention.}

\noindent
$\bullet$
 {\em Everywhere in the present work
  convergence and all relations between random
variables are meant with probability one w.r.t. the measure
$P$ unless specified otherwise. \\
{$\bullet$} A sequence of random
variables $(\zeta_t)_{t\ge1}$ has a property  {\underline {{\bf \em eventually}}} if for
every $\omega$ in a set $\Omega_0$ of $P$ probability 1, the realisation 
 $\zeta_t(\omega)$  has this property for all $t$ greater than some
$t_0(\omega)<\infty$.}\\
{$\bullet$} {\em All conditional expectations are calculated  as integrals w.r.t. corresponding regular conditional probability measures.}\\
{$\bullet$}
{\em The $\inf_{z\in U} h(z)$ of a real valued function $h(z)$ is  $1$
whenever $U=\emptyset$.}

%

%%%%%%%%%%%%%%%%%%%%%%%%%%%%%%%%%%%%%%%%%%%%%%%%%%%%%%%%%%%%%%%%%%%%%%%%%%
                 %       SUB S E C T I O N    Main   Results
%%%%%%%%%%%%%%%%%%%%%%%%%%%%%%%%%%%%%%%%%%%%%%%%%%%%%%%%%%%%%%%%%%%%%%%%%%

 \subsection{Notes on convergence}\label{MR}
%%%%%%%%%%%%%%%%%%%%%%%%%%%%%%%%%%%%%
% Classical stochastic approximation with moving bounds

%%%%%%%%%%%%%%%%%%%%%%%%%%%%%%%%%
  \begin{remark}\label{PoD}{\rm
  This subsection contains  simple results describing sufficient conditions for convergence and rate of convergence. We decided to present this material here  for the sake of completeness, noting that the proof, as well as 
 a number of different sets of sufficient conditions, can be found  in Sharia (2014\nocite{Shar4}) and Sharia and Zhong (2016\nocite{Sh-Zh1}).  
 } \end{remark}
\begin{proposition}\label{SC}   Suppose that  $Z_t$  is a  process  defined by 
\eqref{TSA}, $U_t$ are admissible truncations for $z^0$.

\smallskip
\noindent
$\bullet$  Suppose that 
\begin{description} 
 \item[(D1)] for large $t$'s 
 $$
 (z-z^0)^T R_t(z) \le 0 \;\;\;\mbox{if}\;\;\;z \in U_{t-1};
 $$
\item[(D2)]there exists a predictable process $r_t>0$ such that
$$
  \sup_{z \in U_{t-1}}\frac {E \left\{ \|R_t(z)+\ve_t(z)\|^2 \mid{\cf}_{t-1}\right\}}
  {1+\| z-z^o\|^2}\le r_t $$    
 eventually, and
 $$
  \sum_{t=1}^{\infty}   {r_{t}}{a_t^{-2}} <\infty, \qquad P\mbox{-a.s.} .
$$
 \end{description}
Then $\|Z_t-z^0\|$ converges  ($P$-a.s.) to a finite limit.

\smallskip
\noindent 
$\bullet$ Furthermore, if 
\begin{description}
 \item[(D3)]   for each $\epsilon\in (0, 1),$ there exists a predictable process $\nu_t>0$ such that
$$
 \inf_{\stackrel{ \epsilon \le \|z-z^o\| \le 1/\epsilon}{z\in U_{t-1}}}  -(z-z^0)^T R_t(z)> \nu_t
$$
eventually, where
$$
  \sum_{t=1}^{\infty}  {\nu_{t}}{a_t^{-1}} =\infty, \qquad P\mbox{-a.s.}
$$
\end{description}
Then $Z_t$ converges  ($P$-a.s.)  to $z^0$.

\smallskip
\noindent 
$\bullet$ 
Finally, if  
\begin{description}\item[(W1)] 
$$
\Delta_{t-1}^T R_t(Z_{t-1}) \leq -{\frac 1 2} \Delta a_t\|\Delta_{t-1}\|^2
$$ 
eventually;
\item[(W2)] there exist $0<\delta\leq1$ such that,
$$
\sum_{t=1}^\infty a_t^{\delta-2}E\left\{\|(R_t(Z_{t-1})+{\ve}_t(Z_{t-1}))\|^2 \mid{\cf}_{t-1}\right\}<\infty.
$$
\end{description}
Then $a_t^{\delta}\|Z_t-z^0\|^2$ converges to a finite limit ($P$-a.s.).
\end{proposition}
{\bf Proof.} See Remark \ref{PoD} above.
%%%%%%%%%%%%%%%%%%%%%%%%%%%%%%%
%%%%%%%%%%%%%%%%%%%%%%%%%%%%%%%%%%%%%%%%%%%%%%%%%%%
%
%Asymptotic linearity
%
%%%%%%%%%%%%%%%%%%%%%%%%%%%%%%%%%%%%%%%%%%%%%%%%%%%%

\subsection{Asymptotic linearity }\label{ALin}

In this subsection we establish that under certain conditions, the SA process defined by \eqref{TSA} is asymptotically linear in the statistical sense, that is, it can be represented as a weighted sum of random variables. Therefore, a suitable form of the central limit theorem can be applied to derive  the corresponding asymptotic distribution. 
\begin{theorem}\label{ASM}
Suppose that process $Z_t$ is defined by \eqref{TSA} and 
\begin{description}
\item[(E1)]
\begin{equation}\label{NoTrun}
Z_t=Z_{t-1}+\gamma_t(Z_{t-1})[R_t(Z_{t-1})+\ve_t (Z_{t-1})]\;\;\; \mbox{eventually. }
\end{equation}
\end{description}
Suppose also that there exists a sequence of invertible random matrices $A_t$ such that
\begin{description}
\item[(E2)]
$$
A_t^{-1}\longrightarrow 0
\;\;\;\mbox{ and }\;\;\; 
A_t\gamma_t(z^0) A_t \longrightarrow \eta \;\;\;\mbox{ in probability, }
$$
where $\eta<\infty$ ($P$-a.s.) is a finite matrix;
\item[(E3)]
$$
\lim_{t\rightarrow \infty}A_t^{-1}\sum_{s=1}^t \left [ \Delta \gamma_s^{-1}(z^0) \Delta_{s-1}+\tilde R_s(z^0+\Delta_{s-1}) \right ]=0
$$

in probability, where $$\Delta\gamma_s^{-1}(z^0)=\gamma_s^{-1}(z^0)-\gamma_{s-1}^{-1}(z^0),$$ $$\Delta_s=Z_s-z^0 \;\;\mbox{ and }\;\; \tilde R_s(z)=\gm_s^{-1}(z^0)\gm_s(z)R_s(z);$$
\item[(E4)]
$$
\lim_{t\to \infty}A_t^{-1}\sum_{s=1}^t \Big[\tilde \ve_s(z^0+\Delta_{s-1})-\ve_s(z^0)\Big]=0
$$
in probability, where
$$
\tilde \ve_s(z)=\gm_s^{-1}(z^0)\gm_s(z)\ve_s(z).
$$
\end{description}
Then $A_t(Z_t-Z_t^*)\longrightarrow 0$ in probability where $$Z_t^*=z^0+\gamma_t(z^0) \sum_{s=1}^t\ve_s(z^0);$$ 
that is, $Z_t$ is locally asymptotically linear in $z^0$ with $\gamma_t=\gamma_t(z^0)$ and $\psi_t=\ve_t(z^0)$.
\end{theorem}
{\bf Proof.}
Using the notation $\gamma_t=\gamma_t(z^0)$, $\ve_t=\ve_t(z^0)$ and $\Delta_t=Z_t-z^0$, \eqref{NoTrun} can be rewritten as
$$
\Delta_t-\Delta_{t-1}=\gamma_t \tilde R_t(Z_{t-1})+\gamma_t \tilde\ve_t(Z_{t-1})
$$
eventually. Multiplying both sides by $\gamma_t^{-1}$, we have
$$
\sum_{s=1}^t[\gamma_s^{-1}\Delta_s-\gamma_{s-1}^{-1}\Delta_{s-1}]=\sum_{s=1}^t[\Delta\gamma_s^{-1}\Delta_{s-1}+\tilde R_s(Z_{s-1})+\tilde\ve_s(Z_{s-1})],
$$
and since the sum on the left hand side reduces to $\gamma_t^{-1}\Delta_t-\gamma_0^{-1}\Delta_0$, we obtain
$$
\Delta_t= \gamma_t \left [{\cal H}_t+\sum_{s=1}^t\tilde\ve_s(Z_{s-1}) +\gamma_{0}^{-1}\Delta_{0}\right ]
$$
eventually,
where
$$
{\cal H}_t=\sum_{s=1}^t[\Delta\gamma_s^{-1}\Delta_{s-1}+\tilde R_s(Z_{s-1})].
$$
Since $Z_t-Z_t^*=\Delta_t-(Z_t^*-z^0)$, we have
$$
Z_t-Z_t^*=\gamma_t \Big[{\cal H}_t+\gamma_{0}^{-1}\Delta_{0} \Big]+\gamma_t\sum_{s=1}^t \Big[\tilde \ve_s(Z_{t-1})-\ve_s\Big],
$$
and
$$
A_t(Z_t-Z_t^*)= A_t\gamma_t A_t A_t^{-1} \Big[{\cal H}_t+\gamma_{0}^{-1}\Delta_{0}   \Big ]+A_t\gamma_t A_t A_t^{-1} \sum_{s=1}^t \Big[\tilde \ve_s(Z_{t-1})-\ve_s\Big]
$$
eventually. 
By conditions (E2), (E3) and (E4), we have 
$$
A_t\gamma_t A_t \xrightarrow{P} \eta,\;\;  A_t^{-1} \Big[{\cal H}_t+\gamma_{0}^{-1}\Delta_{0}   \Big ]\xrightarrow{P} 0\;\;\mbox{ and }\;\;A_t^{-1} \sum_{s=1}^t \Big[\tilde \ve_s(Z_{t-1})-\ve_s\Big]\xrightarrow{P} 0
$$
Therefore, $A_t(Z_t-Z_t^*)\longrightarrow 0$ in probability, that is, $Z_t$ is locally asymptotically linear at $z^0$.\hfill $\blacksquare$

\begin{proposition}\label{ASM2}
Suppose that $A_t$ in Theorem \ref{ASM} are positive definite diagonal matrices with non-decreasing elements and
\begin{description}
\item[(Q1)]
$$
A_t^{-2} \sum_{s=1}^t A_s\left [ \Delta\gamma_s^{-1}(z^0)\Delta_{s-1}+\tilde R_s(z^0+\Delta_{s-1})\right ]\longrightarrow 0
$$
\end{description}
in probability, where $\tilde R_t$ is defined in (E3). Then (E3) in Theorem \ref{ASM} holds.
\\
\end{proposition}
{\bf Proof.}
Denote
$$
\chi_s=A_s[\Delta \gamma_s^{-1}(z^0)\Delta_{s-1}+\tilde R_s(z^0+\Delta_{s-1})]
$$
and
$$
A_t^{-1}\sum_{s=1}^t [\Delta \gamma_s^{-1}(z^0)\Delta_{s-1}+\tilde R_s(z^0+\Delta_{s-1})]=A_t^{-1}\sum_{s=1}^t A_s^{-1}\chi_s\;.
$$
Let us denote
$P_s=A_s^{-1}$ and $Q_s=\sum_{m=1}^s\chi_m$. Then using the formula (summation by parts)
$$
\sum_{s=1}^tP_s\Delta Q_s=P_t Q_t-\sum_{s=1}^t\Delta P_s Q_{s-1}\;\;\; \mbox{with}\;\;\;Q_0=0\;,
$$
we obtain
$$
A_t^{-1}\sum_{s=1}^t A_s^{-1}\chi_s=A_t^{-2} \sum_{s=1}^t \chi_s +{\cal G}_t
~~~~ \mbox{where} ~~~~
{\cal G}_t= -A_t^{-1} \sum_{s=1}^t \Delta A_s^{-1}\sum_{m=1}^{s-1}\chi_m.
$$
Since $A_s$ are diagonal, 
$$\Delta A_s^{-1}=A_s^{-1}-A_{s-1}^{-1}=-A_s^{-1}(A_s-A_{s-1})A_{s-1}^{-1}=-\Delta A_s A_s^{-1}A_{s-1}^{-1}.$$ 
Therefore,
$$
{\cal G}_t=A_t^{-1} \sum_{s=1}^t \Delta A_s \left \{A_s^{-1}A_{s-1}^{-1}\sum_{m=1}^{s-1}\chi_m\right\}\;.
$$
Denote by $A_s^{(j,j)}$ the $j$-th diagonal element of $A_s$. Since $0\leq A_{s-1}^{(j,j)}\leq A_{s}^{(j,j)}$ for all $j$,
$$
A_{s-1}^{-2}\sum_{m=1}^{s-1}\chi_m\longrightarrow 0\implies A_s^{-1}A_{s-1}^{-1}\sum_{m=1}^{s-1}\chi_m\longrightarrow 0.
$$
Because of the diagonality, we can apply the Toeplitz Lemma to the elements of ${\cal G}_t$, which gives
$$
A_t^{-1}\sum_{s=1}^t [\Delta \gamma_s^{-1}(z^0)\Delta{s-1}+\tilde R_s(z^0+\Delta_{s-1})]=A_t^{-2} \sum_{s=1}^t \chi_s +{\cal G}_t\longrightarrow 0\;.
$$
\hfill$\blacksquare$
%Remark: center limit theory for i.i.d. and martingale.
\begin{proposition}\label{ASM3}
Suppose that $A_t$ in Theorem \ref{ASM} are positive definite diagonal matrices with non-decreasing elements. Denote by $\alpha^{(j)}$ the $j$-th element of $\alpha\in \mathbb R^m$ and by $A^{(j,j)}$ the $j$-th diagonal element of matrix $A$. Suppose also that 
\begin{description}
\item[(Q2)]
$$
\lim_{t\to\infty} (A_t^{(j,j)})^{-2} \sum_{s=1}^t E\Big\{\Big[\tilde\ve_s^{(j)}(z^0+\Delta_{s-1})-\ve_{s}^{(j)}(z^0)\Big]^2\Big|{\cal F}_{s-1}\Big\}= 0
$$
\end{description}
in probability $P$ for all $j=1,...,m$, where $\tilde\ve_s$ is defined in (E4). Then (E4) in Theorem \ref{ASM} holds.
\end{proposition}
{\bf Proof.} 
Denote $M_t=\sum_{s=1}^t \Big[\tilde \ve_s(z^0+\Delta_{s-1})-\ve_s(z^0)\Big]$. By the assumptions, $M_t$ is a martingale and the quadratic characteristic $\langle M^{(j)}\rangle_t$ of the $j$th component $M_t^{(j)}$ is
$$
\langle M^{(j)}\rangle_t=\sum_{s=1}^t E_{z^0}\Big\{\Big[\tilde\ve_s^{(j)}(z^0+\Delta_{s-1})-\ve_{s}^{(j)}(z^0)\Big]^2\Big|{\cal F}_{s-1}\Big\}.
$$
Using the Lenglart-Rebolledo inequality (see e.g., Liptser and Shiryayev (1989\nocite{LipShir}), Section 1.9), we have
$$
P\Big\{(M_t^{(j)})^2\geq K^2(A_t^{(j,j)})^2\Big\}\leq \frac{\epsilon}{K}+P\Big\{\langle M^{(j)}\rangle_t\;\;\geq\epsilon(A_t^{(j,j)})^2\Big\}
$$
for each $K>0$ and $\epsilon>0$. Now by (Q2), $\langle M^{(j)}\rangle_t/(A_t^{(j,j)})^2\longrightarrow 0$ in probability $P$ and therefore $M_t^{(j)}/A_t^{(j,j)}\longrightarrow0$ in probability $P$. Since $A_t$ is diagonal, (E4) holds.\hfill$\blacksquare$

\begin{remark}\label{ChoNor} {\rm
Let us use Condition (E3) in Theorem \ref{ASM} to construct an optimal step-size sequence $\gm_t(z^0)$. Consider condition (Q1) in the one-dimensional case. Since $R_t(z^0)=0$, we have
\begin{eqnarray*}
&&A_t\left[\Delta \gamma_t^{-1}(z^0)\Delta_{t-1}+\tilde R_t(z^0+\Delta_{t-1})\right]\\
&=&\left[\Delta \gamma_t^{-1}(z^0)+e_t\frac{ R_t(z^0+\Delta_{t-1})- R_t(z^0)}{\Delta_{t-1}}\right]A_t\Delta_{t-1},
\end{eqnarray*}
where $e_t=\gm_t^{-1}(z^0)\gm_t(z^0+\Delta_{t-1})$. In most applications, the rate of $A_t$ is $\sqrt t$ and $\sqrt t \Delta_{t}$ is stochastically bounded. Therefore, for (Q1) to hold, one should at least have the convergence 
$$
\Delta \gamma_t^{-1}(z^0)+e_t\frac{ R_t(z^0+\Delta_{t-1})- R_t(z^0)}{\Delta_{t-1}}\longrightarrow 0.
$$
If $\gm_t(z)$ is continuous, given that $\Delta_t\longrightarrow0$, we expect $e_t\longrightarrow1$. Therefore, we should have
$$
\Delta \gamma_t^{-1}(z^0)\approx -R_t'(z^0).
$$
Using the similar arguments for the multi-dimensional cases, we expect the above relation to hold for large $t$'s, where $R_t'(z^0)$ is the matrix of the derivatives of $R_t(z)$ at $z=z^0$. So, an optimal choice of the step-size sequence should be
$$
\gm_t^{-1}(z)=-\sum_{s=1}^t R_s'(z),
$$
or a sequence which is asymptotically equivalent to this sum.
} 
\end{remark}
%%%%%%%%%%%%%%%%%%%%%%%%%%%%%%%%%%%
%
% ev
%%%%%%%%%%%%%%%%%%%%%%%%%%%%%%
\begin{remark}\label{ev}{\rm 
{\bf (a)} Condition (E1) in Theorem \ref{ASM} holds  if the truncations in \eqref{TSA} do not occur for large $t$'s. More precisely, (E1) holds if  for $t > T$  the truncations in \eqref{TSA}  do not occur for some, possibly random $T$. 

\noindent
{\bf (b)}
Let us now consider the case when $U_t$ is a shrinking sequence. For example, suppose that a consistent, but not necessarily efficient, auxiliary estimator $\tilde Z_t$ is available. Then one can take the truncations on $U_t=S(\tilde Z_t, r_t)$, which is a sequence of closed spherical sets in $\mathbb R^m$ with the center at $\tilde Z_t$ and the radius $r_t\longrightarrow0$. The resulting  procedure is obviously consistent, as $\| Z_t-\tilde Z_t\|\leq r_t\longrightarrow 0$ and $\tilde Z_t\longrightarrow  z^0$.  However, if  $r_t$ decreases too rapidly, condition (E1)
may fail to hold. Intuitively,  it is quite obvious that we should not allow $r_t$ to decreases too rapidly, as it may result in $ Z_t$  having  the same asymptotic  properties as $\tilde Z_t$, which might not be optimal. This truncation will be admissible if $\|\tilde Z_t- z^0\|<r_t$ eventually. In these circumstances, (E1) will hold if the procedure generates the sequence  $ Z_t$ which converges to $z^0$ faster than $r_t$ converges to 0.

 \noindent
{\bf (c)} The considerations described  in (b) lead to the following construction. Suppose that an auxiliary estimator $\tilde Z_t$ has a convergence rate $d_t$, in the sense that $d_t$ is a sequence of positive r.v.'s such that $d_t\longrightarrow \infty$ and $d_t(\tilde Z_t- z^0) \to 0$ ~  $P$-a.s. Let us consider the following truncation sets
$$
U_t=S\left(\tilde Z_t,c(d_t^{-1}+a_t^{-1})\right),
$$
where $c$ and $a_t$ are positive and $a_t\longrightarrow\infty$.  Then the truncation sequence is obviously admissible since $\|\tilde Z_t- z^0\| < cd_t^{-1}$ eventually. Now, if  we can claim (using Proposition \ref{SC} or otherwise)  that $a_t\| Z_t - z^0\| \longrightarrow 0 $,  then condition (E1) holds. Indeed, suppose that (E1) does not hold, that is, 
the truncations in \eqref{TSA} occur infinitely many times on a set $A$ of positive probability. This would imply that $ Z_t$ appears on the surface of the spheres $U_t$  infinitely many times on $A$. Since $ z^0 \in S(\tilde Z_t, c d_t^{-1})$ eventually, we obtain that $\| Z_t- z^0\| \geq c a_t^{-1}$ infinitely many times on $A$, which contradicts our assumptions.

Another possible choice of the truncation sequence is
$$
U_t=S\left(\tilde Z_t,c\left(d_t^{-1} \vee a_t^{-1}\right)\right).
$$
(Here, $a\vee b=\max(a,b)$ and $a\wedge b=\min(a,b)$). 
If  we can claim by Proposition \ref{SC} or otherwise  that
$a_t\| Z_t - z^0\| \to 0 $,  then condition 
(E1) holds.  Indeed, suppose that (E1) does not hold, that is,     on a set $A$ of positive probability
the truncations in \eqref{TSA} occur infinitely many times. This would imply that
$$
\|\tilde Z_t- Z_t\|=c(d_t^{-1} \vee a_t^{-1})
$$
and
$$
1= c^{-1}(d_t \wedge a_t)  \|\tilde Z_t- Z_t\| \le
c^{-1}(d_t \wedge a_t) \|\tilde Z_t- z^0\|+ c^{-1}(d_t \wedge a_t)\| Z_t- z^0\|
$$
infinitely many times on $A$, which contradicts our assumptions.
}\end{remark}

%%%%%%%%%%%%%%%%%%%%%%%%%%%%%%%%%%%%
%
%  Special models and examples
%
%%%%%%%%%%%%%%%%%%%%%%%%%%%%%%%%%%%
\section{Special models and examples}\label{SpME}

\subsection{Classical  problem of stochastic approximation}\label{CSA}

Consider the classical problem of stochastic approximation to find a root $z^0$ of the equation $R(z^0)=0$. 
 Note that in the classical case, the step-size sequence can in general be of the form form $\gamma_t (Z_{t-1})=a_t ^{-1}  \gamma(Z_{t-1})$. 
 However, without loss of generality we can assume that  $\gamma_t= a_t^{-1}\bf I$, since  $\gamma(Z_{t-1})$ can be included in $R$ and $\ve_t$.
Therefore, taking the step-size sequence $\gamma_t= a_t^{-1}\bf I$, where $a_t\longrightarrow \infty$ is a predictable scalar process, let us  consider the procedure 
\begin{equation}\label{SN}
Z_{t}=\Phi_{U_t}\Big(Z_{t-1}+ a_t^{-1} [ R(Z_{t-1})+\ve_t(Z_{t-1})]\Big).
\end{equation}
  \begin{remark}\label{PoD} {\rm
 In the corollary below we derive simple sufficient conditions for asymptotic linearity in the case when  $a_t=t $.  We also assume, using Proposition \ref{SC} or otherwise,  that
 $t^{\delta/2} (Z_t-z^0) \longrightarrow 0$ for any $\delta \in (0,1)$.   Note also that the condition (A1) below requires that the procedure is designed in such a way  that the truncations in \eqref{SN} do not occur for large $t$'s (see Remark \ref{ev} for a detailed discussion of this requirement).}
 \end{remark}

\begin{corollary} \label{CAS}
Suppose that $Z_t$ is defined by \eqref{SN}, $a_t=t$ and $t^{\delta/2} (Z_t-z^0) \longrightarrow 0$ for any $\delta \in (0,1)$. Suppose also that
\begin{description}
\item[(A1)]
$$
Z_{t}=Z_{t-1}+\frac1 {t} [ R(Z_{t-1})+\ve_t(Z_{t-1})] \;\;\;\mbox{ eventually;}
$$
\item[(A2)]
$$
R(z^0+u)=-u+\alpha(u)
\;\;\;\mbox{ where }\;\;\;
\|\alpha (u)\|= O(u^{1+\epsilon})
$$
as $u \to 0$ for some $\epsilon>0$;
\item[(A3)]
$$
t^{-1}\sum_{s=1}^tE\Big\{\Big[\ve_s(z^0+u_s)-\ve_s(z^0)\Big]^2\Big|{\cal F}_{s-1}\Big\}<\infty,
$$
where $u_s$ is any predictable process with the property $u_s\longrightarrow 0$.
\end{description}
Then $Z_t$ is asymptotically linear.
\end{corollary}
%%%%%%%%%%%%%
{\bf Proof.}
Let $A_t=\sqrt t {\bf I}$, then $A_t \gamma_t A_t={\bf I}$ since $\gm_t={\bf I}/t$. Condition (E2) in Theorem \ref{ASM} is satisfied. On the other hand, since $\tilde R(z)=R(z)$ and $\Delta\gamma_t^{-1}={\bf I}$, we have
$$
A_t^{-2} \sum_{s=1}^t A_s\left [ \Delta\gamma_s^{-1}\Delta_{s-1}+\tilde R_s(Z_{s-1})\right ]=\frac 1 t \sum_{s=1}^t \sqrt s [\Delta_{s-1}+R(z^0+\Delta_{s-1})]=\frac 1 t \sum_{s=1}^t \sqrt s \alpha(\Delta_{s-1}).
$$
By (A2), there exists a constant $K>0$ such that 
$$
\|\sqrt s \alpha (\Delta_{s-1})\|\leq K\left\|\sqrt s\Delta_{s-1}^{1+\epsilon}\right\|
=K\left\|\sqrt{\frac s {s-1}}\left [(s-1)^{\frac{1} {2(1+\epsilon)}}\Delta_{s-1}\right]^{1+\epsilon}\right\|
$$
eventually.
Since $1/ [2(1+\epsilon)]<1/2$, we have $(s-1)^{1/ {[2(1+\epsilon)]}}\Delta_{s-1}\longrightarrow 0$, and therefore $\|\sqrt s \alpha (\Delta_{s-1})\|\longrightarrow 0$ as $\Delta_s\longrightarrow 0$. Thus, by the Toeplitz Lemma (see Lemma \ref{Toep} in Appendix A),
$$
\frac 1 t \sum_{s=1}^t \sqrt s \alpha(\Delta_{s-1})\longrightarrow 0
$$
So, (Q2) in Proposition \ref{ASM2} holds implying that condition (E3) in Theorem \ref{ASM} is satisfied. Since $\tilde \ve _t(z)=\ve _t(z)$, it follows from (A3) that condition (Q2) in Proposition \ref{ASM3} holds. This implies that (E4) in Theorem \ref{ASM} holds. Thus, all the conditions of Theorem \ref{ASM} hold, implying that $Z_t$ is asymptotically linear. \hfill$\blacksquare$

%%%%%%%%%%%%%

\begin{remark} {\rm Using asymptotic linearity, the asymptotic normality is an immediate consequence of Corollary \ref{CAS}. Indeed, we have $\sqrt t(Z_t-Z_t^*)\longrightarrow 0$ in probability, where
$$
Z_t^*=z^0+\frac 1 t\sum_{s=1}^t \ve_s(z^0).
$$

\noindent So, $Z_t$ and $Z_t^*$ have the same asymptotic distribution. Now, to obtain the asymptotic distribution of $Z_t$, it remains only to apply the central limit theorem for martingales.
}\end{remark}
%%%%%%%%%%%%%%%%%%%%%%%%
\begin{remark} {\rm Note that condition (A2) above  assumes that $R$ function should be scaled in such a way that the derivative at $z^0$ is $-1$.
Alternatively,   a step-size sequence should be considered of the form  $\gamma_t (Z_{t-1})=t ^{-1} \gamma(Z_{t-1}),$ with appropriately chosen 
$ \gamma(Z_{t-1}).$ Detailed discussion of selection of an appropriate step-size sequence in the context of statistical parametric estimation is given in Section \ref{PEGM}.
}
\end{remark}
%%%%%%%%%%%%%%%%%%%%%%%%%%%%%% 
%
%          Finding a root of a polynomial
% 
%%%%%%%%%%%%%%%%%%%%%%%%%%%%%%%%%%%%%%%%%%%%%%%%%%%%%%%%%%%%%%%

\begin{example}\label{Poly}
{\rm 
Let  $l$  be a positive integer and 
 $$
 R(z)=-\sum_{i=1}^{l}C_i(z-z^0)^i,
 $$                
where $z , z^0 \in \mathbb{R}$ and $C_i$ are real constants. Suppose that
$$
(z-z^0)R(z)\leq0 \;\;\;\mbox{ for all } \;\;\;z\in \mathbb R.
$$
Unless  $l =1$, we cannot use the standard  SA   without truncations as the standard condition on the rate of growth at infinity  does not hold.  So, we consider $Z_t$ defined by \eqref{SN} with a  slowly 
expanding truncation sequence $U_t=[-u_t, u_t]$, where 
 $$
    \sum_{t=1}^{\infty}  u_t^{2l}~ a_t^{-2} <\infty.
$$
We can assume for example, that   $u_t=Ct^{r/{2l}}$, where $C$ and $r$ are some positive constants and $r < 1$.  One can also take a truncation sequence which is independent of $l$, e.g.,  $u_t=C \log t$, where $C$ is a positive constant.

Suppose for simplicity that  the measurement errors are state free with the property that 
$ \sum_{t=1}^{\infty} \sigma_t^2 {a_t^{-2}} <\infty $, where 
$
\sigma_t^2= {E \left\{ \ve_t^2 \mid{\cf}_{t-1}\right\}}.
$
Then $|Z_t-z^0|$ converges ($P$-a.s.) to a finite limit. Furthermore, if $z^0$ is a unique root, 
then  $Z_t\longrightarrow z^0$ ($P$-a.s.) 
provided that 
$
\sum_{t=1}^{\infty}   a_{t}^{-1}=\infty.
$
  Finally, if  $Z_t$ is defined by \eqref{SN} with  $a_t=C_1t$, then $t^\alpha (Z_t-z^0)\xrightarrow{a.s.} 0$ for any $\alpha<1/2$ (see  Sharia and Zhong (2016\nocite{Sh-Zh1}) for details). So, it follows that conditions in Corollary \ref{CAS} hold (with  $R$ replaced by $C_1^{-1}R$),  implying that   $Z_t$ is locally asymptotically linear. Now, depending on the nature of the error terms, one can apply a suitable form of the central limit theorem to obtain  asymptotic normality of $Z_t$.
}
\end{example}

%%%%%%%%%%%%%%%%%%%%%%%%%%%%%%
%
%          Linear procedures
% 
%%%%%%%%%%%%%%%%%%%%%%%%%%%%%%%%%%%%%%%%%%%%%%%%%%%%%%%%%%%%%%%

\subsection{Linear procedures}
Consider the recursive procedure
\begin{equation}\label{LP}
Z_t=Z_{t-1}+\gamma_t(h_t-\beta_t Z_{t-1})
\end{equation}
where $\gamma_t$ is a predictable positive definite matrix process, $\beta_t$ is a predictable positive semi-definite matrix process and $h_t$ is an adapted vector process (i.e., $h_t$ is ${\cal F}_t$-measurable for $t\geq 1$). If we assume that $E\{h_t|{\cal F}_{t-1}\}=\beta_t z^0$, we can view \eqref{LP} as a SA procedure designed to find the common root $z^0$ of the linear functions
$$
R_t(u)=E\{h_t-\beta_t u|{\cal F}_{t-1}\}=E\{h_t|{\cal F}_{t-1}\}-\beta_t u=\beta_t(z^0-u)
$$ 
which is observed with the random noise 
$$
\ve_t=\ve_t(u)=h_t-\beta_t u-R_t(u)=h_t-E\{h_t|{\cal F}_{t-1}\}=h_t-\beta_t z^0.
$$
%%%%%%%%%%%%%%%%
 \begin{remark}\label{PoD}    {\rm Recursive procedures \eqref{LP} are linear in the sense that they locate the common root $z^0$ of the linear functions $R_t(u)=\beta_t(z^0-u)$.  The second part  of the corollary  below shows that the process $Z_t$ is asymptotically linear in the statistical sense, that is, it can be represented as a weighted sum of random variables.  The first part of the corollary  below  contains  sufficient conditions for convergence and rate of convergence. We decided to present this material here  for the sake of completeness, noting that the proof can be found  in  Sharia and Zhong (2016\nocite{Sh-Zh1}) (note also that (G1) below will hold if, e.g.,  $\Delta\gm_t^{-1}=\beta_t$).}
    \end{remark}
    %%%%%%%%%%%%%%%%%%%  corollary 

\begin{corollary}\label{LRA}
Suppose that $Z_t$ is defined by \eqref{LP} with $E(h_t|{\cal F}_{t-1})=\beta_tz^0$ and   $a_t$ is a non-decreasing positive predictable process. 

\smallskip
\noindent
{\large  \bf  \em 1.}  Suppose that 
\begin{description}
\item[(G1)]
$\Delta\gamma_t^{-1}-2\beta_t+\beta_t\gamma_t\beta_t$ is negative semi-definite eventually;
\item[(G2)] 
$$
\sum_{t=1}^{\infty} a_t^{-1}E\{(h_t-\beta_t z^0)^T\gamma_t(h_t-\beta_t z^0)|{\cal F}_{t-1}\}<\infty.
$$
\end{description}
Then $a_t^{-1}(Z_t-z^0)^T\gamma_t^{-1}(Z_t-z^0)$ converges to a finite limit (P-a.s.). 

\smallskip
\noindent
{\large \bf \em 2.}
Suppose that  $\gamma_t\longrightarrow 0$ and
\begin{equation}\label{CLP}
\gamma_t^{1/2}\sum_{s=1}^t (\Delta\gamma_s^{-1}-\beta_s)\Delta_{s-1}\longrightarrow 0
\end{equation}
in probability, where $\Delta_t=Z_t-z^0$. 

Then $Z_t$ is asymptotically linear, that is,
$$
\gm_t^{1/2}(Z_t-z^0)=\gm_t^{-1/2}\sum_{s=1}^t \ve_s+r_t(z^0),
$$
where $r_t(z^0)\longrightarrow 0$ in probability. 
\\
\end{corollary}

%%%%%%%%%%%%%%%%%%%%%%%%%%%%%%%
\noindent
{\bf Proof.}
Let us check the conditions of Theorem \ref{ASM} for $A_t=\gamma_t^{-1/2}.$ Conditions (E1) and (E2) trivially hold. Since $\ve_t(u)=h_t-\beta_tz^0$ is state free (i.e. does not depend on $u$), (E4) also holds. Since
$
\tilde R_t(Z_{t-1})=R_t(Z_{t-1})
=-\beta_t\Delta_{t-1},
$
we have 
$$
A_t^{-1}\sum_{s=1}^t \left ( \Delta \gamma_s^{-1}(z^0) \Delta_{s-1}+\tilde R_s(z^0+\Delta_{s-1}) \right )=\gamma_t^{1/2}\sum_{s=1}^t (\Delta\gamma_s^{-1}-\beta_s)\Delta_{s-1}\;,
$$
and (E3) now follows from \eqref{CLP}. Thus, all conditions of Theorem \ref{ASM} are satisfied which implies the required result.\hfill $\blacksquare$

 %%%%%%%%%%%%%%%%%%%%%%%%%%%%%%%%%%%%%%
 %
%                                     AR1  example
%
%%%%%%%%%%%%%%%%%%%%%%%%%%%%%%%%%%%%%%

 \begin{example}\label{AR1}
{\rm   Corollary \ref{LRA}   can be applied to study asymptotic behaviour of recursive least squares  estimators in regression or time series models.
To demonstrate this, let us consider  a simple example of AR(1) process
$$
X_t=\theta X_{t-1}+\xi_t,
$$
where
${\xi}_t$ is a sequence of square integrable random variables with mean zero.  
Consider the  recursive  least squares (LS)  estimator of $\theta$ defined by
%(ArLsq)
                          \begin{eqnarray}
&&\hat\theta_t=\hat\theta_{t-1}+\hat I_t^{-1}
X_{t-1}\left(X_t-\hat \theta_{t-1}X_{t-1}\right) ,   \nonumber \\
&&\hat I_t=\hat I_{t-1}+X_{t-1}^2, \qquad t=1,2,\dots \nonumber
\end{eqnarray}
where $\hat\theta_0$  and  $\hat I_0 > 0$ are any starting points and 
$
\hat I_t=\hat I_0+\sum_{s=1}^t X_{s-1}^2.
$
This procedure is clearly a particular case of \eqref{LP} with
$$
 z^0=\theta, ~~~~~ Z_t= \hat\theta_t, ~~~~~  \gamma_t=\hat I_t^{-1}, ~~~~~h_t= X_{t-1}X_t, ~~~~~~~ \beta_t=X_{t-1}^2.
$$
 Since 
 $\Delta\gm_t^{-1}=X_{t-1}^2= \beta_t$, condition (G1)  holds (see Corollary 5.2 in Sharia and Zhong (2016\nocite{Sh-Zh1})). Also,  since  
 $$
 h_t-\beta_t z^0=X_{t-1}(X_t -X_{t-1}\theta) =X_{t-1} \xi_t,
 $$
 it follows that 
  $$
E\{(h_t-\beta_t z^0)^T\gamma_t(h_t-\beta_t z^0)|{\cal F}_{t-1}\}= X_{t-1}^2 \hat I_t^{-1} E\{\xi_t^2 |{\cal F}_{t-1}\}.
$$ 
 Let $0<\delta<1$. Then taking $a_t=\hat I_t ^\delta$ in (G2) we obtain 
  $$
\sum_{t=1}^{\infty} a_t^{-1}E\{(h_t-\beta_t z^0)^T\gamma_t(h_t-\beta_t z^0)|{\cal F}_{t-1}\}=
 \sum_{t=1}^{\infty} \frac 1{\hat I_t^{1+\delta}} X_{t-1}^2 E\{\xi_t^2 |{\cal F}_{t-1}\}
$$ 
Now, since $\Delta\hat I_t=X_{t-1}^2$, if   $ \hat I_t \to \infty$ then the sum above is finite even if the conditional variances  
 $E\{\xi_t^2 |{\cal F}_{t-1}\}$  go to infinity with rate  $ \hat I_t^{\delta^0}$, as far as  ${\delta^0 <\delta}$ (this trivially follows from, e.g., Lemma 6.3 in 
 Sharia and Zhong (2016\nocite{Sh-Zh1})). 
 
 Let us now assume for simplicity that ${\xi}_t$ is a sequence of i.i.d. r.v.'s with mean zero and variance $1$.  Then 
 consistency  and rate of convergence follows 
 without any further moment assumptions  on the innovation process. Indeed,  since  $ \hat I_t \to \infty$ for any $\theta\in \mathbb{R}$
 (see, e.g, Shiryayev (1984\nocite{Shir}, Ch.VII, $\S$5),  it follows that 
  all the conditions of  part 1 in Corollary \ref{LRA} hold implying that $I_t^{1+\delta} (\hat\theta_t-\theta)^2$ converges a.s. to a finite limit for any  
  $0<\delta<1$ and $\theta\in\mathbb{R}$. 
  
Furthermore,  since
 $\Delta\gm_t^{-1}= \beta_t$, \eqref{CLP} trivially holds. It  therefore follows that  $\hat\theta_t$ is asymptotically linear and  asymptotic normality is
 now obtained by applying the central limit theorem  for i.i.d. random variables. 
 }
\end{example}

%%%%%%%%%%%%%%%%%%%  corollary 

%%%%%%%%%%%%%%%%%%%%%%%%%%%%%%%%%%%%%%%%%%%%%
%                                             Application to parameter estimation 
%
%%%%%%%%%%%%%%%%%%%%%%%%%%%%%%%%%

\subsection{Application to parameter estimation}\label{PEGM}
Let $X_1, \dots, X_n$ be  random variables with a
joint distribution depending on an unknown parameter $\theta$. Then an
$M$-estimator of $\theta$ is defined as  a solution of the
estimating equation
%%%%%%%%%%%%%%%%%%%%%%%%%%%%%%%%%%%%%%%%%%%%%%%%%%%%
%                  \eqref{esteqg}
%%%%%%%%%%%%%%%%%%%%%%%%%%%%%%%%%%%%%%%%%%%%%%%%%%%%%%
\begin{equation}\label{esteqg}
\sum_{i=1}^n \psi_i(\theta)=0,
\end{equation}
where   $\psi_i(\theta)=\psi_i(X_1^i; \theta),$  ~$i=1,2,\dots, n$,
are suitably chosen functions which may, in general, depend on the  vector
$X_1^i=(X_1, \dots,X_i )$ of all past and present observations.
  If $f_i(x,\theta)=f_i(x,\theta|X_1, \dots,X_{i-1})$ is the
conditional probability density function or probability function
of the observation $X_i,$ given $X_1, \dots,X_{i-1},$   then one
can obtain a   MLE (maximum likelihood estimator)  on choosing
\begin{equation}\label{CML}
\psi_i(\theta)=l_t(\theta)= [f_i'(\theta, X_i|X_1^{i-1})]^T/f_i(\theta, X_i|X_1^{i-1}).
\end{equation} 
Besides MLEs, the class of
$M$-estimators includes estimators
 with special properties such as {robustness}.
Under certain regularity and ergodicity conditions, it can be
proved that there exists a consistent sequence of solutions of
\eqref{esteqg} which has  the property of local asymptotic
linearity.

  Let us consider
 estimation procedures which are recursive in the sense that each successive
 estimator is obtained from the previous one by a simple adjustment.
   In particular, we consider  a class of
estimators
%%%%%%%%%%%%%%%%%%%%%%%%%%%%%%%%%%%%%%%%%%%%%%%%%%
%                   \eqref{rec1}
%%%%%%%%%%%%%%%%%%%%%%%%%%%%%%%%%%%%%%%%%%%%%%%%%
$$
\hat\theta_t=\Phi_{U_t}\Big[\hat\theta_{t-1}+
{\gamma_t(\hat\theta_{t-1})}\psi_t(\hat\theta_{t-1})\Big],
~~~~~~~~~t\ge 1,
$$
where $\psi_t$ is a suitably chosen vector process, $\gamma_t$  is
a   matrix valued step-size process, and
$\hat\theta_0\in {\mathbb{R}}^m$ is an initial value.
This type of recursive estimators are especially convenient when the corresponding $\psi$-functions are non-linear in $\theta$  and therefore, solving \eqref{esteqg}  would require a numerical method (see e.g., Example \ref{Gamex}).  A detailed discussion and a heuristic justification
of this estimation procedure are given  in Sharia (2008)\nocite{Shar1}.

The above procedure can be rewritten in the SA form. Indeed, assume that $\theta$ is an arbitrary but fixed value of the parameter and denote
$$
R_t(z)=E_{\theta}\left\{ \psi_t(z)\mid{\cf}_{t-1}\right\}
~~~ \mbox{and} ~~~ 
\ve_t(z)=\left(\psi_t(z) - R_t(z)\right).
$$  
Following the argument in Remark \ref{ChoNor} (see also Sharia (2010\nocite{Shar3})), an optimal step-size sequence would be
$$
\gamma_t^{-1}(\theta)=-\sum_{s=1}^t R'_s(\theta)
$$
 If
$\p_t(z)$ is differentiable w.r.t. $z$ and differentiation
of $ R_t(z)=E_{\theta} \{ \p_t(z)\mid {\cf}_{t-1}\}$
is allowed under the integral sign, then $R'_t(z)=E_{\theta} \{ \p_t'(z)\mid {\cf}_{t-1}\}.$
 This implies that,
for a given sequence of estimating functions  $\psi_t(\theta),$ another possible
 choice of the step-size sequence is
$$
\gm_t(\theta)^{-1}=-\sum_{s=1}^t E_{\theta} \{ \p_s'(\theta)\mid
{\cf}_{s-1}\},
$$
or any sequence with the increments
$$
\Dl \gm_t^{-1}(\theta)=\gm_t^{-1}(\theta)-\gm_{t-1}^{-1}(\theta)= -E_{\theta} \{ \p_t'(\theta)\mid
{\cf}_{t-1}\}.
$$
Also, since $\psi_t(\theta)$ is typically a $P^{\theta}$-martingale difference, 
$$
%E_{\theta}\{ \p_t(\theta)\mid {\cf}_{t-1}\}
0=\int \p_t(\theta,x\mid X_1^{t-1}) f_t(\theta,x\mid X_1^{t-1})\mu
(dx),
$$
and if the differentiation w.r.t. $\theta$ 
 is allowed under the integral sign, then (see Sharia (2010\nocite{Shar3}) for details)
$$ 
 E_{\theta}\{ \p_t'(\theta)\mid {\cf}_{t-1}\}
 =-E_{\theta}\{ \p_t(\theta)l^{T}_t(\theta)\mid
 {\cf}_{t-1}\},
$$
where  $l_t(\theta)$ is defined in \eqref{CML}.
Therefore, another possible choice of the step-soze sequence is any sequence with the increments
$$
\Dl \gm_t^{-1}(\theta)=\gm_t^{-1}(\theta)-\gm_{t-1}^{-1}(\theta)=E_{\theta}\{ \p_t(\theta)l^{T}_t(\theta)\mid
 {\cf}_{t-1}\}.
$$
Therefore, since the process
$$
M_t^\theta=\sum_{s=1}^t\p_s(\theta)
$$
is a $P^\theta$-martingale, the above sequence can be rewritten as
$$
\gm_t^{-1}(\theta)=\langle M^\theta, U^\theta \rangle_t
$$
where $U_t^\theta=\sum_{s=1}^t l_s(\theta)$  is the  score martingale.

Let us consider a likelihood case,
 that is   $\psi_t(\theta)=l_t(\theta)$, the above sequence
 is the conditional Fisher
information 
\begin{equation}
I_t(\theta)=\sum_{s=1}^t E\{l_s(\theta)l_s^T(\theta)|{\cal F}_{s-1}\}.
\end{equation}
Therefore, the corresponding recursive procedure is
%%%%%%%%%%%%%%%%%%%%%%%%%%%%%%%%%%%%%%%%%%%%%%%%%%
%                   \eqref{recmle}
%%%%%%%%%%%%%%%%%%%%%%%%%%%%%%%%%%%%%%%%%%%%%%%%%
\begin{equation} \label{recmle}
\hat\theta_t=\Phi_{U_t}\Big(\hat\theta_{t-1}+
{I_t^{-1}(\hat\theta_{t-1})}l_t(\hat\theta_{t-1})\Big),
~~~~~~~~~t\ge 1,
\end{equation}
Also, given that the model possesses certain  ergodicity properties,
asymptotic linearity of \eqref{recmle} implies asymptotic
efficiency.   In particular, in the case
of i.i.d. observations, it follows that the above recursive procedure is
asymptotically
normal with parameters $(0, \; i^{-1}(\theta) )$, where $i(\theta)$ is the one-step Fisher information.

\subsubsection{The i.i.d case}
Consider the classical scheme of i.i.d. observations $X_1, X_2,...$ having a common probability density function $f(x,\theta)$ w.r.t. some $\sigma$- finite measure $\mu$, where $\theta\in\mathbb R^m$. Suppose that $\psi(x,\theta)$ is an estimating function with
$$
E_\theta\left\{\psi(X_1,\theta)\right\}=\int \psi(x,\theta)f(x,\theta)\mu(dx)=0.
$$ 
A recursive estimator $\hat{\theta}_t$ can be defined by
$$
\hat\theta_t=\Phi_{U_t}\Big(\hat\theta_{t-1}+a_t^{-1}
{\gm(\hat\theta_{t-1})}\psi(X_t,\hat\theta_{t-1})\Big)
$$
where $a_t$ is a non-decreasing real sequence, $\gm(\theta)$ is an invertible $m\times m$ matrix and truncation sequence $U_t$ is admissible for $\theta$. 
In most applications $a_t=t$ and an optimal choice of $\gamma(\theta)$ is
$$
\gamma(\theta)=\Big[E_\theta\Big\{\psi(X_t,\theta)l^T(X_t,\theta)\Big\}\Big]^{-1} ~~ \mbox{where}~~
l(x,\theta)=\frac{[f'(x,\theta)]^T}{f(x,\theta)}\;.
$$

%%%%%%%%%%%%%%%%%%%%%%%%%%%%%%%%%%%%%%%%
%example of Gamma distribution  \label{Gamex}
%%%%%%%%%%%%%%%%%%%%%%%%

\begin{example}\label{Gamex}
Let  $X_1,X_2,\ldots$ be i.i.d. random variables from Gamma$(\theta,1)$ ($\theta >0$).
Then the  the common probability density function  is
$$
f(x,\theta)=\frac1{ {\bf{\Gamma}(\theta)} } x^{\theta-1}e^{-x}, \;\;\;
\theta >0, \;\; x >0,
$$
where ${\bf {\Gamma}(\theta)}$ is the Gamma function. Denote  
$$
{{\log}' {\bf{\Gamma}(\theta)}}=\frac{d}{d\theta}{\log} {\bf{\Gamma}(\theta)},\;\;\;{{\log}'' {\bf{\Gamma}(\theta)}}=\frac{d^2}{d\theta^2}{{\log}} {\bf{\Gamma}(\theta)}.
$$
Then 
$$
\frac {f'(x,\theta)}{f(x,\theta)}={\log} x-
{\log}' {\bf{\Gamma}(\theta)}
  ~~~~~\mbox{ and }~~~~
  {i} (\theta)={{\log}''} {\bf{\Gamma}(\theta)},
$$
where $ {i} (\theta)$ is the one-step Fisher  information.
Then a recursive likelihood estimation procedure can be defined as
%                 EstG
 \begin{equation} \label{EstG}                      
\hat \theta_t=\Phi_{U_t}\left(\hat \theta_{t-1}+\frac{1}{t ~ {\log}'' {\boldsymbol{\Gamma}(\hat\theta_{t-1})}}\left[
\log X_t-{{\log}' {\boldsymbol{\Gamma}(\hat\theta_{t-1})}}\right]
\right)
\end{equation}
with $U_t=[\alpha_t,\beta_t]$ where
$\alpha_t\downarrow  0$ and $\beta_t\uparrow \infty $ are sequences of positive numbers. Then it can be shown that (see Appendix B) if
 \begin{equation} \label{albt}         \sum_{t=1}^\infty
\frac {\alpha_{t-1}^2} {t}= \infty   ~~~~~ \mbox{and} ~~~~~ \sum_{t=1}^{\infty} \frac{\log^2\alpha_{t-1}
+\log^2\beta_{t-1}}{t^2}  < \infty,
\end{equation}
then $\hat\theta_t$ is strongly consistent and asymptotically efficient, i.e., $
\hat\theta_t\xrightarrow{a.s.} \theta $ as $ t\longrightarrow
\infty, 
$ 
and 
$$
{\cal L}\Big(t^{1/2}(\hat \theta_t-\theta)|P^\theta\Big)\xrightarrow{w}{\cal N}\Big(0,\log''\boldsymbol\Gamma(\theta)\Big).
$$
For instance,  
$$
{\alpha_t=C_1 ({{\log}} \; (t+2))^{-\frac12}} \;\;\; \mbox{and} \;\;\;
{\beta_t=C_2(t+2)}
$$
with  some positive constants $C_1$ and $C_2$, obviously satisfy  \eqref{albt}.

The above result can be derived by rewriting   \eqref{EstG}   in the form of the stochastic approximation (see Appendix B for details), i.e.,
%                 SapG
 \begin{equation} \label{SapG}                      
\hat \theta_t=\Phi_{U_t}\left(\hat \theta_{t-1}+\frac1{t} \left[ R(\hat\theta_{t-1}) +\ve_t(\hat\theta_{t-1}) \right]
\right)
\end{equation}
where 
$$
R(u)=R^\theta(u)=\frac{1}{{\log}'' {\boldsymbol{\Gamma}(u)}}E_\theta\{log X_t-\log' {\bf {\Gamma}}(u)\}=
\frac{1}{{\log}'' {\boldsymbol{\Gamma}(u)}}\left(\log' {\bf {\Gamma}}(\theta)-
\log' {\bf {\Gamma}}(u)\right)
$$
and 
$$
\ve_t(u)=\frac{1}{{\log}'' {\boldsymbol{\Gamma}(u)}}\left[
\log X_t-{{\log}' {\boldsymbol{\Gamma}(u)}}\right] - R(u).
$$
\end{example}

%%%%%%%%%%%%%%%%%%%%%%
%{Simulations}\label{SIMU}
%%%%%%%%%%%%%%%%%%%%%

\section{Simulations}\label{SIMU}
\subsection{Finding roots of polynomials}\label{MCPoly}

Let us consider a problem described in Section \ref{Poly} with
$$
R(z)=-(z-z^0)^7+2(z-z^0)^6-5(z-z^0)^5-3(z-z^0),
$$
and suppose that the random errors are independent Student random variables with degrees of freedom 7. Consider SA procedure \eqref{SN} with $a_t=3t$ and the truncation sequence $U_t=[-\log3t,\log3t]$. Then (see Example \ref{Poly}), it follow  that this procedure is consistent, i.e., converges almost surely to $z^0$, and asymptotically linear.  Also, since the error terms are i.i.d., it follows that the procedure is asymptotically normal. Note that the SA without truncations fails to satisfy the standard condition on the rate of growth at infinity. Here, slowly expanding truncations are used to artificially slow down the growth of $R$ at infinity.

Figure 1 shows 30 steps  of the procedure with starting points at $-2$, $0$ and $5$ respectively, where the root $z^0=2$.
 A histogram of the estimator  over 500 replications (with $Z_0=0$) is shown in Figure 2.
\begin{figure}[h!]
  \label{PolyMov}
   \centering
   \includegraphics[width=0.7\textwidth]{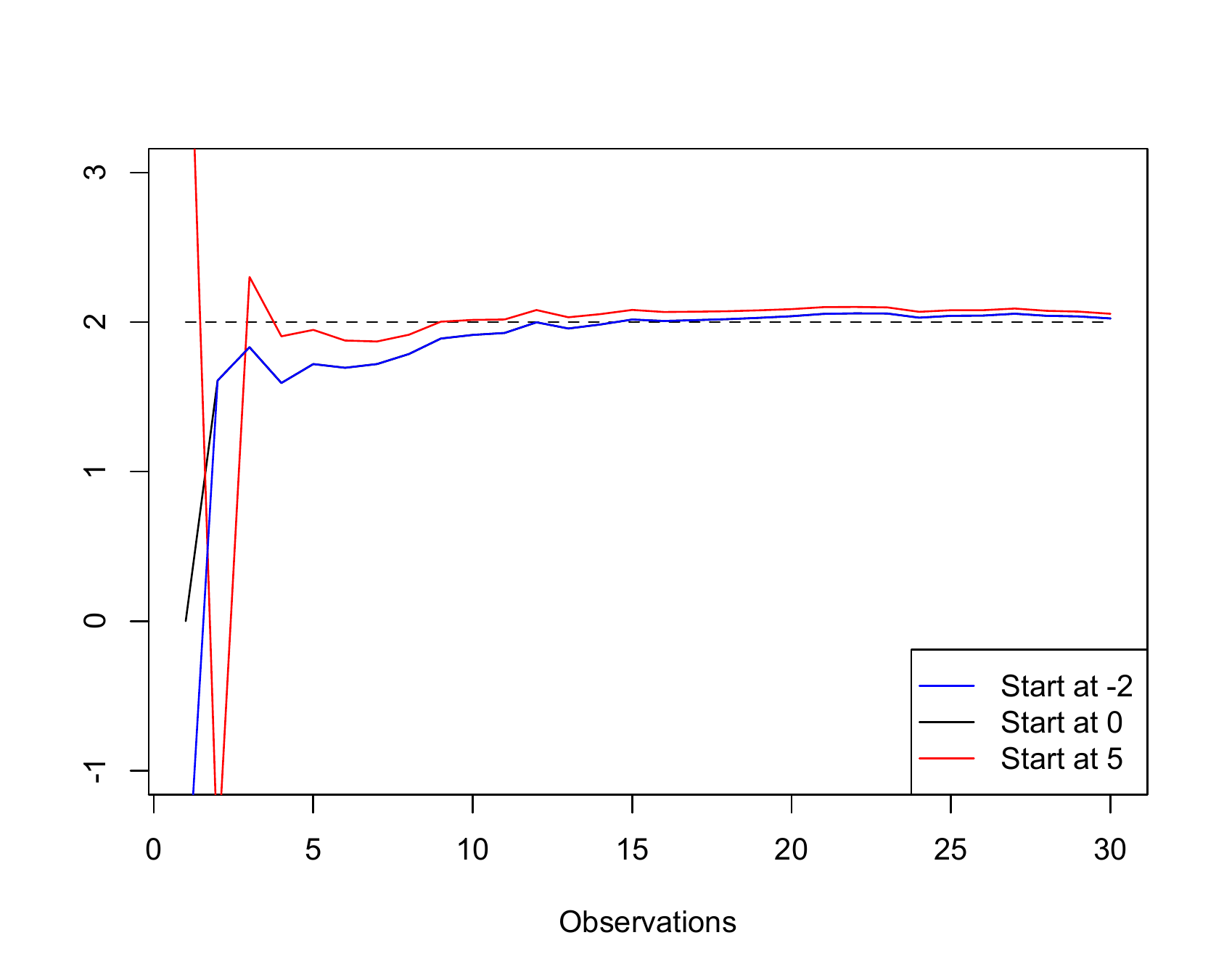}
  \caption{Realizations of the estimator in the polynomial example}
\end{figure}
\begin{figure}[h!]
  \label{PolyHist}
   \centering
   \includegraphics[width=0.6\textwidth]{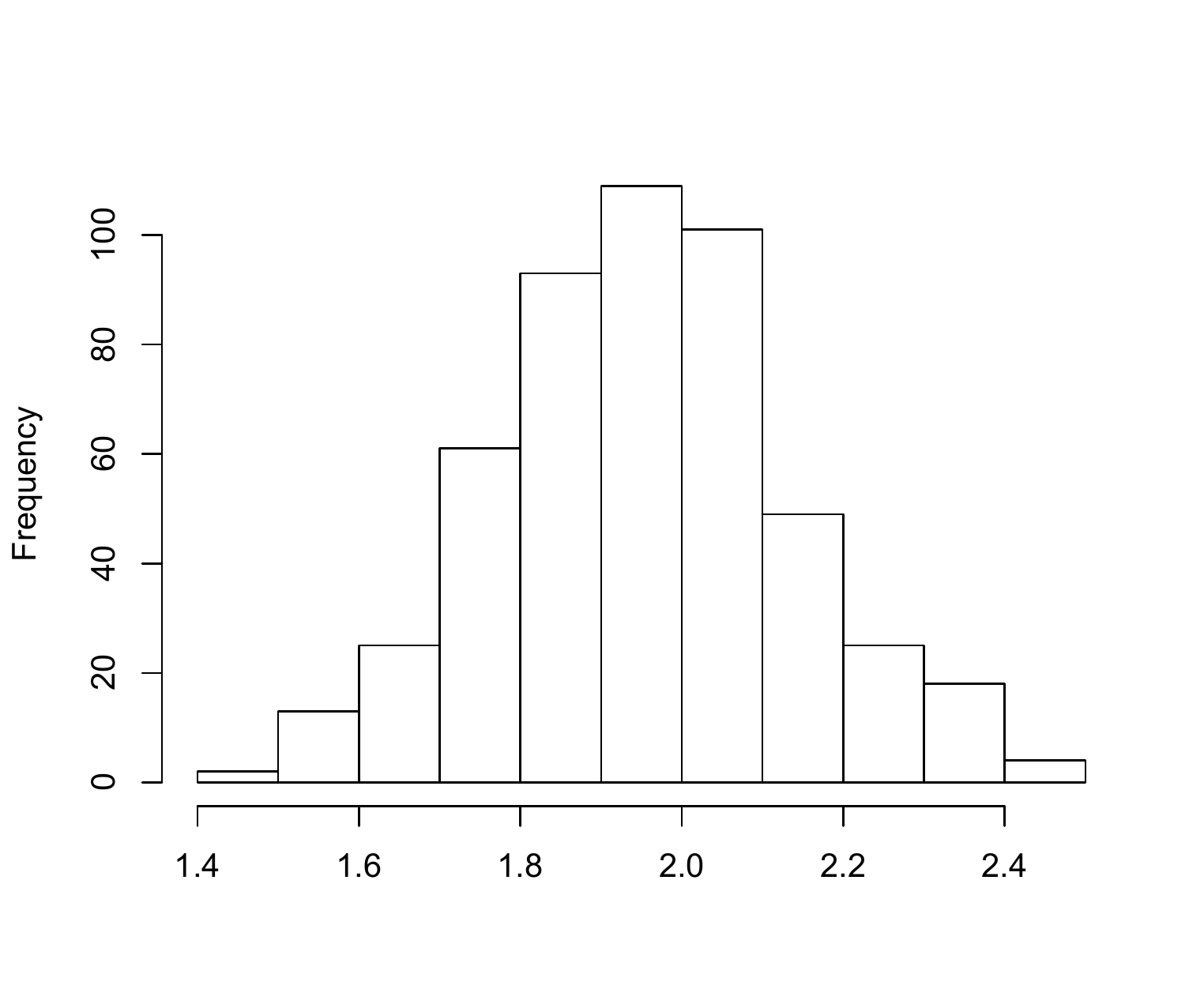}
  \caption{Histogram of the estimator in the polynomial example}
\end{figure}
\subsection{Estimation of the shape parameter of the Gamma distribution}
Let us consider procedure \eqref{EstG} in Example \ref{Gamex} with following two sets of truncations $U_t=[\alpha_t,\beta_t]$.
\begin{description}
\item[(1)] FT -- Fixed truncations: $\alpha_t=\alpha$ and $\beta_t=\beta$ where $0<\alpha<\beta<\infty$.
\item[(2)] MT -- Moving truncations: $\alpha_t=C_1[\log(t+2)]^{(-1/2)}$ and $\beta_t=C_2(t+2)$ where $C_1$ and $C_2$ are positive constants.
\end{description}
\begin{figure}[h!]
  \label{GamNoBoost}
   \centering
   \includegraphics[width=0.7\textwidth]{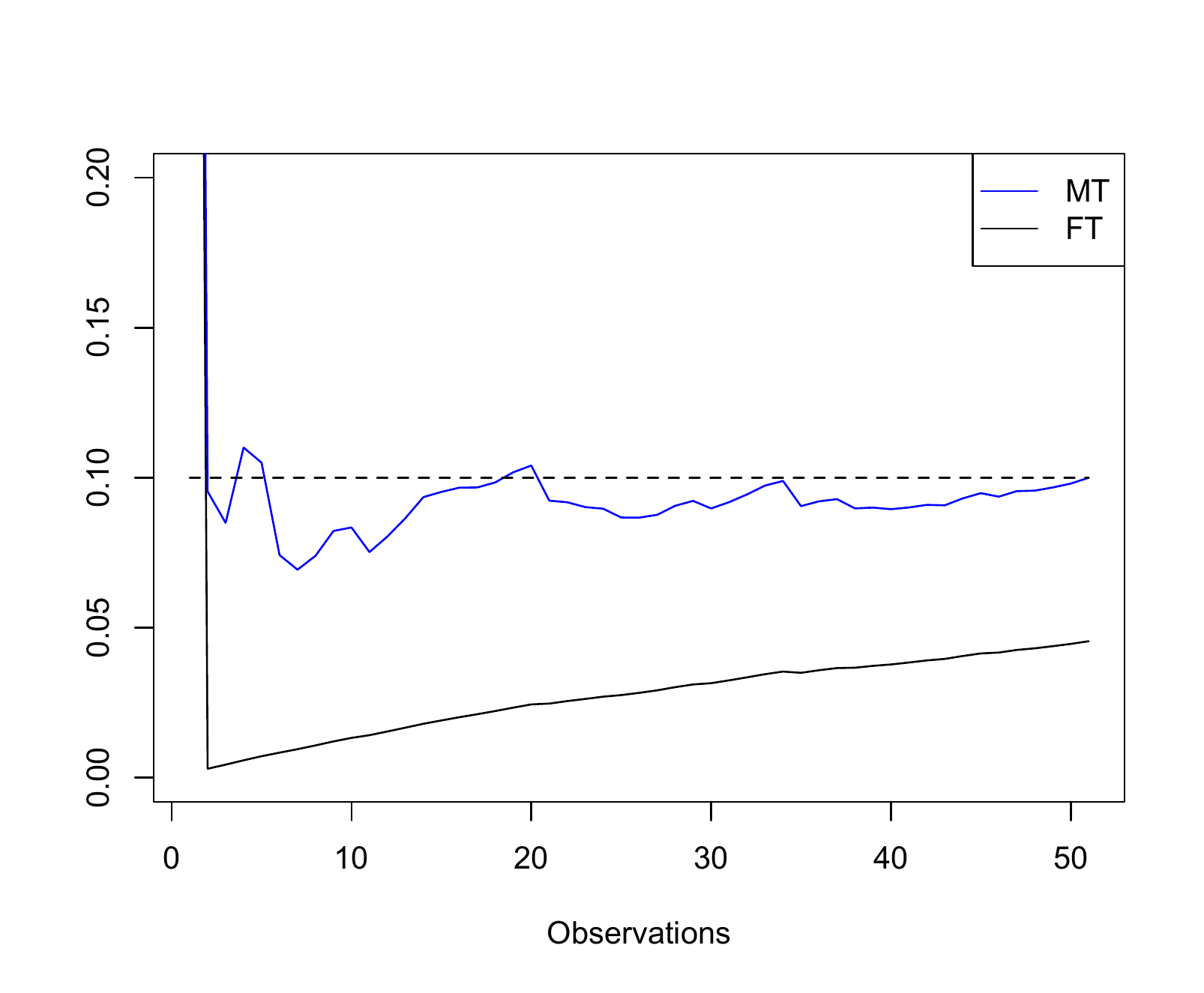}
  \caption{Performance of the estimator of the parameter in the Gamma distribution}
\end{figure}
Figure 3 shows realizations of procedures \eqref{EstG} when $\theta=0.1$ and the starting point $\hat \theta_0=1$, $C_1=0.1$, $C_2=1$ in MT, and $\alpha=0.003$, $\beta=100$ in FT. As we can see, the MT estimator approaches the true value of $\theta$ following a zigzag path. However, the FT estimator moves very slowly towards the true value of $\theta$, caused by singularity at 0 of the functions appearing in the procedure.

%%%%%%
%        {Appendix A}\label{App}
%%%%%%

%\newpage

\numberwithin{equation}{section}
\numberwithin{lemma}{section}
\addcontentsline{toc}{chapter}{Appendix}
\section{Appendix }\label{App}

 \begin{lemma} (The Toeplitz Lemma)\label{Toep}
Let $\{a_n\}$ be a sequence of non-negative numbers such that $\sum_{n=1}^\infty a_n$ diverges. If $\nu_n \longrightarrow \nu_\infty$ as $n \longrightarrow \infty$, then
$$
\lim_{n \longrightarrow \infty} {\frac{\sum_{i=1}^n a_i \nu_i}{\sum_{i=1}^n a_i}} = \nu_\infty\;.
$$
\end{lemma}
{\bf Proof.} Proof can be found in Lo{\`e}ve (1977\nocite{Loeve}, P.250).
\hfill $\blacksquare$
\\
%%%%%%%%%%%%%%%%%%%%%%%%%%%%%%%%%%%%%%%%%%%%%%%%%%%%

\bigskip
\noindent 
{\bf  Properties of Gamma distribution} ~~
In Example \ref{Gamex},
we will need the following properties of  the Gamma function  (see, e.g.,  Whittaker (1927)\nocite{Wit}, 12.16).
${\log}' {\boldsymbol{\Gamma}}$ is  increasing, 
    ${\log}'' {\boldsymbol{\Gamma}}$  is decreasing and continuous, 
$$
{\log}'' {\boldsymbol{\Gamma}}(x) \leq \frac {1+x}{x^2}
$$
and
\begin{equation}\label{Log''G2}
{\log}'' {\boldsymbol{\Gamma}}(x) \geq \frac 1 x.
\end{equation}
Also (see Cramer (1946)\nocite{Cram},  12.5.4), 
$$
{\log}' {\boldsymbol{\Gamma}}(x) \le {\mbox{ln}} (x).
$$
Then,
\begin{equation}\label{expect}
E_\theta\left\{\log  X_1 \right\}= {\log}' {\boldsymbol{\Gamma}}(\theta),            ~~~~ ~~~~
E_\theta\left\{ \left(\log X_1\right)^2 \right\}={\log}'' {\boldsymbol{\Gamma}}(\theta) +
 \left({\log}' {\boldsymbol{\Gamma}}(\theta)\right)^2,
  \end{equation}
$$
E_\theta\left\{ \left(\log X_1 -{\log}' {\boldsymbol{\Gamma}}(\theta)\right)^2 \right\}={\log}'' {\boldsymbol{\Gamma}}(\theta). 
$$
Using  \eqref{Log''G2} and \eqref{expect}  we obtain
\begin{equation}\label{Sq}
E_\theta \left\{ \|R(u)+\ve_t(u)\|^2 \mid{\cf}_{t-1}\right\}=\frac{ {\log}'' {\boldsymbol{\Gamma}(\theta)}
+\left(\log' {\bf {\Gamma}}(\theta)-
\log' {\bf {\Gamma}}(u)\right)^2}   {({\log}'' {\boldsymbol{\Gamma}(u))^2}  }\;.
\end{equation}
\\
The convergence to $\theta$ of the estimator defined by  \eqref{EstG}
is shown in Sharia (\nocite{Shar4}2014). To establish the rate of convergence, let us show that the conditions of Corollary 4.5 in Sharia and Zhong (2016)\nocite{Sh-Zh1} hold. 
Since 
\begin{eqnarray*}
R'(u)=\frac{d R(u)}{d u}&=&-\frac{{\log}'' {\boldsymbol{\Gamma}(u)}}{{\log}'' {\boldsymbol{\Gamma}(u)}}-\frac{{\log}''' {\boldsymbol{\Gamma}(u)}}{[{\log}'' {\boldsymbol{\Gamma}(u)}]^2}\left(\log' {\bf {\Gamma}}(\theta)-
\log' {\bf {\Gamma}}(u)\right)\\
&=&-1-\frac{{\log}''' {\boldsymbol{\Gamma}(u)}}{[{\log}'' {\boldsymbol{\Gamma}(u)}]^2}\left(\log' {\bf {\Gamma}}(\theta)-
\log' {\bf {\Gamma}}(u)\right),
\end{eqnarray*}
we have $R'(\theta)=-1\leq-1/2$ and condition (B1) of Corollary 4.5 in Sharia and Zhong (2016)\nocite{Sh-Zh1} 
%\ref{ClaRC}
 holds. Since $E_{\theta}\left\{ \ve_t(u)\mid {\cf}_{t-1} \right\}=0$, we have 
\begin{equation}\label{axali}
E_{\theta}\left\{ [ R(u)+\ve(u) ]^2 \mid {\cf}_{t-1} \right\} =
R^2(u) +E_{\theta}\left\{  \ve_t^2(u)  \mid {\cf}_{t-1} \right\}. 
\end{equation}
Using \eqref{Sq} and \eqref{axali}, 
\begin{eqnarray*}
E_\theta\left\{  \ve_t^2(u)  \mid {\cf}_{t-1} \right\}&\leq& E_\theta\left\{ [ R(u)+\ve(u) ]^2 \mid {\cf}_{t-1} \right\}\\
&=& {\log}'' {\boldsymbol{\Gamma}(\theta)}
+\left(\log' {\bf {\Gamma}}(\theta)-
\log' {\bf {\Gamma}}(u)\right)^2,
\end{eqnarray*}
which is obviously a continuous function of $u$. Thus, for any $v_t\longrightarrow0$, we have $E_{\theta}\left\{  \ve_t^2(\theta+v_t)  \mid {\cf}_{t-1} \right\}$ converges to a finite limit and so condition (BB) in Corollary 4.7 in Sharia and Zhong (2016)\nocite{Sh-Zh1}  holds. Therefore, all the conditions of this  corollary are satisfied with $a_t=t$ implying that $t^\delta (\hat{\theta}_t-\theta)^2\xrightarrow{a.s.}0$ for any $\delta<1$.

Furthermore, since the second derivative of $R(u)$ exists, $R'(\theta)=-1$, and $R(\theta)=0$, by the Taylor expansion,
$$
R(\theta+u)=-u+R''(\tilde u)u^2
$$
for small $u$'s and for some $\tilde u>0$. Therefore, condition (A2) in Corollary \ref{CAS} holds. It is also easy to check that
$$
E_{\theta}\Big\{\Big[\ve_s(\theta+u_s)-\ve_s(\theta)\Big]^2\Big|{\cal F}_{s-1}\Big\}\longrightarrow 0
$$
for any predictable process $u_s\longrightarrow 0$. Condition (A3) is immediate from the Toeplitz Lemma.
Thus, estimator $\hat{\theta}_t$ defined by \eqref{SapG} is asymptotic linear. Now, using the CLT for i.i.d. r.v.'s, it follows that  $\hat{\theta}_t$ is asymptotically efficient.

%\nocite{Shar0}\nocite{Shar}\nocite{Shar5}

%\nocite{AnM}\nocite{Andr}\nocite{AMP}\nocite{Bena}\nocite{Ben}\nocite{Bor}\nocite{BaS}\nocite{Chen1}\nocite{Chen2}\nocite{Cam}\nocite{Cram}\nocite{Dov}\nocite{Eng}\nocite{EaH}\nocite{Fab}\nocite{Gu}\nocite{Hay}\nocite{HaJ}\nocite{Kall}\nocite{Khas}\nocite{KaW}\nocite{Kush1}\nocite{Kush2}\nocite{Lai}\nocite{Laz}\nocite{Lel}\nocite{Lju}\nocite{LaS}\nocite{MaM}\nocite{PaJ}\nocite{Pol}\nocite{RM}\nocite{Rob2}\nocite{Sac}\nocite{Sakr}\nocite{Shar0}\nocite{Shar}\nocite{Shar1}\nocite{Shar2}\nocite{Shar3}\nocite{Shar4}\nocite{Shar5}\nocite{Shir}\nocite{Tadic1}\nocite{Tadic2}\nocite{WaD}\nocite{Wit}\nocite{Loeve}\nocite{Blum}\nocite{Glady}\nocite{RobSieg}\nocite{NevKhas}\nocite{Serf}\nocite{Huber}\nocite{LehCas}\nocite{LipShir}\nocite{Feigin}\nocite{kushner3}\nocite{kushnerCla}\nocite{kushnerShw}\nocite{burk}\nocite{hodgesLeh}\nocite{chung}\nocite{laiRob}\nocite{laiRob2}\nocite{wei}\nocite{delyon}\nocite{laiYing}\nocite{laz2}
 
%\newpage

\bibliographystyle{acm}
\bibliography{bibLT1}

\end{document}